\documentstyle[a4,leqno,amsfonts,12pt]{article}

\newcommand{\C}{\mbox{\bf C}}

\newcommand{\rbox}{$\:\:$ \raisebox{-1ex}{$\:\Box\:$}}
\newcommand{\OC}{\overline{\mbox{\bf C}}}

\newcommand{\N}{\mbox{\bf N}}
\newcommand{\bP}{\mbox{\bf P}}

\newcommand{\de}{\delta}
\newcommand{\mb}{\mbox}

\newcommand{\beq}{\begin{equation}}
\newcommand{\eeq}{\end{equation}}
\newcommand{\oge}{\succeq}
\newcommand{\ole}{\preceq}
\newcommand{\ve}{\varepsilon}

\newcommand{\ov}{\overline}
\newcommand{\Om}{\Omega}
\newcommand{\om}{\omega}
\newcommand{\z}{\zeta}
\newcommand{\ka}{\kappa}
\newcommand{\ga}{\gamma}

\newtheorem{th}{Theorem}

\newtheorem{lem}{Lemma}

\newcommand{\ueberschrift}{\bigskip\goodbreak\noindent\bigskip}
\newcounter{theabsatz}
\newcommand{\absatz}[1]{\stepcounter{theabsatz} \ueberschrift
                           {\large \bf \arabic{theabsatz}. {#1}} \setcounter{equation}{0}}

\begin{document}
\mathsurround=2pt

\begin{center}
{\large \bf
SIMULTANEOUS APPROXIMATION AND\\ INTERPOLATION
OF FUNCTIONS\\
 ON CONTINUA IN THE COMPLEX PLANE}\\[3ex]
VLADIMIR V. ANDRIEVSKII\\
IGOR E. PRITSKER\\
RICHARD S. VARGA\\[3ex]
\end{center}
\footnotetext{This research was supported in part by the National Science
Foundation grant DMS-9707359.
}

\begin{abstract}
We construct polynomial approximations of  Dzjadyk type (in terms of the $k$-th
modulus of continuity, $k \ge 1$) for analytic functions defined on a continuum
$E$ in the complex plane, which simultaneously interpolate at given points of
$E$. Furthermore, the error in this approximation is decaying as
$e^{-cn^\alpha}$ strictly inside $E$, where $c$ and $\alpha$ are positive
constants independent of the degree $n$ of the approximating polynomial.
\end{abstract}



\absatz{Introduction and main results}

Let $E\subset\C$ be a compact set with connected complement $\Om:=\OC\setminus
E$, where $\OC:=\C\cup\{\infty\}$ is the extended complex plane. Denote by $A(E)$
the class of all functions continuous on $E$ and analytic in $E^0$, the interior
of $E$ (the case $E^0=\emptyset$ is not excluded). Let $\bP_n,\,
n\in\N_0:=\{ 0,1,2,\ldots\},$ be
the class of complex polynomials of degree at most
$n$. For $f\in A(E)$ and $n\in\N_0$, define
$$
E_n(f,E):=\inf\limits_{p\in\bP_n}||f-p||_E\, ,
$$
where $||\cdot||_E$ denotes the uniform norm on $E$. By Mergelyan's theorem
(see \cite{dzj}), we have that
$$
\lim\limits_{n\to\infty}E_n(f,E)=0\quad (f\in A(E)).
$$
The following
assertion on ``simultaneous approximation and interpolation"
quantifies a result of Walsh \cite[p. 310]{wal}:
Let $z_1,\ldots,z_N\in E$ be distinct points, $f\in A(E)$.
Then for any $n\in\N:=\{ 1,2,\ldots\},\, n\ge N-1$, there exists a polynomial
$p_n\in\bP_n$ such that

\beq
\label{1.1}
||f-p_n||_E\le c\, E_n(f,E),
\eeq
$$
p_n(z_j)=f(z_j)\quad (j=1,\ldots,N),
$$
where $c>0$ is independent of $n$ and $f$.

A suitable polynomial has the form
$$
p_n(z)=p_n^*(z)+\sum_{j=1}^N\frac{q(z)}{q'(z_j)(z-z_j)}
(f(z_j)-p_n^*(z_j)),
$$
where
$$
q(z):=\prod_{j=1}^N(z-z_j),
$$
and $p_n^*\in\bP_n$ satisfies
$$
||f-p_n^*||_E=E_n(f,E).
$$
It is natural to ask whether it is possible to interpolate the function $f$
as before at arbitrary prescribed points and to simultaneously  approximate
it in an even stronger sense than in (\ref{1.1}). The theorem of Gopengauz
\cite{gop} about simultaneous polynomial approximation of real functions continuous
on the interval $[-1,1]$ and their interpolation at $\pm 1$ is an example of
such result. For recent accounts of improvements and generalizations of this
remarkable statement (for real functions) we refer the reader to \cite{pri},
\cite{szaver} and \cite{kilpre}.

We shall make use of the D-approximation (named after Dzjadyk, who found in the
late 50's - early 60's a constructive description of H\"older classes requiring a
nonuniform scale of approximation) as a substitute for (\ref{1.1}). There is an
extensive bibliography devoted to this subject (see, for example, the monographs
\cite{dzj}, \cite{tam}, \cite{gai}, \cite{she} and \cite{andbeldzj}). In the
overwhelming majority of the results on D-approximation, $E$ is a continuum
(one of the rare exceptions is the recent interesting paper \cite{shi1}). In
\cite{and1} it is shown that, for the D-approximation to hold for a continuum
$E$, it is sufficient and under some mild restrictions also necessary that $E$
belongs to the class $H^*$, which can be defined as follows (cf. \cite{and2}
and \cite{and3}).

From now on we assume that $E$ is a continuum with diam$E>0$, connected
complement $\Om$ and boundary $L:=\partial E$. In the sequel, we denote by
$\alpha,\,\beta,\, c,c_1,\ldots$ positive constants (possibly different at
different occurences) that either are absolute or depend on
parameters not essential for the arguments; otherwise, such a dependence will be
indicated.

We say that $E\in H$ if any points $z,\z\in E$ can be joined by an arc $\ga
(z,\z)\subset E$ whose length $|\ga(z,\z)|$ satisfies the condition
\beq
\label{1.3} |\ga(z,\z)| \le c\, |z-\z|,\quad c=c(E)\ge 1.
\eeq
Let us compactify the domain $\Om$ by prime ends in the Caratheodory sense (see
\cite{pom1}). Let $\tilde{\Om}$ be this compactification, and let $\tilde
L:=\tilde{\Om}\setminus\Om$. Assuming that $E\in H$, then all the prime ends
$Z\in\tilde L$ are of the first kind, i.e., they have singleton impressions
$|Z|=z\in L$. The circle $\{\xi:\, |\xi-z|=r\},\, 0<r< \frac{1}{2}$diam$E$,
contains one arc, or finitely many arcs, dividing $\Om$ into two subdomains: an
unbounded subdomain and a bounded subdomain such that $Z$ can be defined by a
chain of cross-cuts of the bounded subdomain. Let $\ga_Z(r)$ denote that one of
these arcs for which the unbounded subdomain is as large as possible (for given
$Z$ and $r$). Thus, the arc $\ga_Z(r)$ separates the prime end $Z$ from
$\infty$ (cf. \cite{bel}, \cite{andbeldzj}).

If $0<r<R<\frac{1}{2} $ diam$E$, then $\ga_Z(r)$ and $\ga_Z(R)$ are the
sides of some quadrilateral $Q_Z(r,R)\subset \Om$ whose other two
sides are parts of the boundary $L$. Let $m_Z(r,R)$ be the module of
this quadrilateral, i.e., the module of the family of arcs that separate the
sides $\ga_Z(r)$ and $\ga_Z(R)$ in $Q_Z(r,R)$ (see \cite{ahl}, \cite{lehvir}).

We say that $E\in H^*$ if $E\in H$ and if there exist constants $c=c(E)<\frac{1}{2}
$diam$E$ and $c_1=c_1(E)$ such that
\beq
\label{1.4}
|m_Z(|z-\z|,c)-m_{\cal Z}(|z-\z|,c)| \le c_1
\eeq
for any pair of prime ends $Z,{\cal Z} \in {\tilde L},$ with their impressions $z=|Z|,\,
\z=|\cal Z|$ satisfying $|z-\z|<c$.

In particular, $H^*$ includes domains with quasiconformal boundaries
(see \cite{ahl}, \cite{lehvir}) and the classes $B_k^*$ of domains introduced by
Dzjadyk \cite{dzj}. For a more detailed investigation of the geometric
meaning of conditions (\ref{1.3}) and (\ref{1.4}), see \cite{and3}.

We will be studying functions defined by their $k$-th modulus of continuity
$(k\in\N)$. There is a number of different definitions
of these moduli in the complex plane (see \cite{vorpol}, \cite{tam},
\cite{dyn}, \cite{she1}). The definition by Dyn'kin \cite{dyn} is the most
convenient for our purpose here.

From now on, suppose that $E\in H^*$. Set
$$ D(z,\de):=\{\z:\, |\z-z|\le\de\}\quad (z\in\C,\, \de>0). $$
The quantity
$$ \om_{f,k,z,E}(\de):=E_{k-1}(f,E\cap D(z,\de)), $$
where $f\in A(E),\, k\in\N,\, z\in E,\, \de>0$, is called the
{\it $k$-th local modulus of continuity}, and
$$ \om_{f,k,E}(\de):=\sup_{z\in E}\om_{f,k,z,E}(\de) $$
is called the {\it $k$-th (global) modulus of continuity} of
$f$  on $E$. It is known (see \cite{tam}) that the behavior of this modulus is
essentially the same as in the classical case of the interval $E=[-1,1]$. In
particular,
\beq \label{1.5}
\om_{f,k,E}(t\de)\le c\, t^k\,\om_{f,k,E}(\de)\quad (t>1,\de>0).
\eeq

We denote by $A^r(E),\, r\in\N$, the class of functions $f\in A(E)$ which are
$r$-times continuously differentiable on $E$, where we set $A^0(E):=A(E)$.

By definition,
the function $w=\Phi(z)$ maps $\Om$ conformally and univalently onto
$\Delta:=\{
w:\, |w|>1\}$ and is normalized by the conditions
\beq
\label{1.55}
\Phi(\infty)=\infty,\,\, \Phi'(\infty)>0.
\eeq
The same symbol $\Phi$ denotes the homeomorphism between the compactification
$\tilde{\Om}$ of $\Om$ and $\ov{\Delta}$, which coincides with $\Phi(z)$ in $\Om$.
Let $\Psi:=\Phi^{-1}$. We define the distance to the level curves of $\Phi(z)$
$$
L_\de:=\{\z:\, |\Phi(\z)|=1+\de\}\quad (\de>0)
$$
by
$$
\rho_\de(z):=\mb{dist}(z,L_\de)\quad (z\in\C,\, \de>0),
$$
where
$$
\mb{dist}(\z,B):=\inf\{|\z-z|: z \in B\}\quad (\z\in\C,\,B\subset\C).
$$

\begin{th}
\label{th1}
Let $E\in H^*,\, f\in A(E),\,  k\in\N$, and let
$z_1,\ldots,z_N\in E$ be distinct points. Then for any $n\in\N,\, n\ge N+k$,
there exists a polynomial $p_n\in\bP_n$ such that
\beq
\label{1.6}
|f(z)-p_n(z)|\le c_1\,
\om_{f,k,E}(\rho_{1/n}(z))\quad (z\in L),
\eeq
\beq
\label{1.7}
p_n(z_j)=f(z_j)\quad (j=1,\ldots,N)
\eeq
with $c_1$ independent of $n$.

Moreover, if $E^0\not=\emptyset$ and if for any $0<\de<1$, there is a constant
$c_2$ such that
\beq
\label{1.77}
\int\limits_0^\de  \om_{f,k,E}(t)\,\frac{dt}{t}\le c_2\,
\om_{f,k,E}(\de),
\eeq
then, in addition to (\ref{1.6}) and (\ref{1.7}),
\beq
\label{1.8}
||f-p_n||_K\le c_3\, \exp(-c_4n^\alpha)
\eeq
for every compact set $K\subset E^0$,
where the  constants $c_3,c_4$ and $0<\alpha\le 1$ are independent of $n$.
\end{th}

A polynomial $p_n$ satisfying (\ref{1.6}) is called a D-approximation of the
function $f$  (D-property of $E$, Dzjadyk type direct theorem). For $k>1$,
(\ref{1.6}) generalizes the corresponding direct theorems of Belyi and Tamrazov
\cite{beltam} (when $E$ is a quasidisk) and Shevchuk  \cite{she1} (when $E$
belongs to the Dzjadyk class $B_k^*$). More detailed history can be found in
these papers.

It was first noticed by Shirokov \cite{shi3} that the rate of D-approximation
may  admit significant improvement strictly inside $E$. Saff and Totik
\cite{saftot} proved that if $L$ is an analytic curve, then an exponential rate
is achievable strictly inside $E$, while on the boundary the approximation is
``near-best". However, even for domains with piecewise smooth boundary without
cusps (and therefore belonging to $H^*$), the error of approximation strictly
inside $E$ cannot be better than $e^{-cn^\alpha}$ (cf. (\ref{1.8})), where
$\alpha$ may  be arbitrarily small (see \cite{mai}, \cite{shitot}). In the
results from \cite{mai}, \cite{shitot}, \cite{shi2} containing estimates of the
form (\ref{1.8}), it is usually assumed that $\Om$ satisfies a wedge condition.
For a continuum $E\in H^*$, this condition can be violated.

Keeping in mind the Gopengauz result \cite{gop}, we generalize Theorem
\ref{th1}  to the case of the Hermite interpolation and simultaneous
approximation of a function $f\in A^r(E)$ and its derivatives. For
simplicity we formulate and prove this assertion only for the case
of boundary interpolation points and without the analog of
(\ref{1.8}).
\begin{th}
\label{th2d}
Let $E\in H^*,\, f\in A^r(E),\, r\in\N,\, k\in\N$, and let
$z_1,\ldots,z_N\in \partial E$
be distinct points. Then for any $n\in\N,\, n\ge Nr+k$,
there exists a polynomial $p_n\in\bP_n$ such that
for $l=0,\ldots,r$,
\beq
\label{1.6d}
|f^{(l)}(z)-p^{(l)}_n(z)|\le c\, \rho_{1/n}^{r-l}(z)\,
\om_{f^{(r)},k,E}(\rho_{1/n}(z))\quad (z\in L),
\eeq
and
\beq
\label{1.7d}
p^{(l)}_n(z_j)=f^{(l)}(z_j)\quad (j=1,\ldots,N),
\eeq
with $c$ independent of $n$.
\end{th}

Our next goal is to allow the number of interpolation nodes $N$
to grow infinitely with the degree of approximating polynomial
$n$. It is well known that we cannot take $N-1$ equal to $n$,
preserving uniform convergence (cf. Faber's theorem \cite{fab}
claiming that for $E=[-1,1]$ there is no universal set of nodes
such that the Lagrange interpolating polynomials converge to
every continuous function in uniform norm). However, it was first
observed by Bernstein \cite{ber} that for any continuous function
on $E=[-1,1]$ and any small $\ve>0$, there exists a sequence of
polynomials interpolating in the Chebyshev nodes and uniformly
convergent on $[-1,1]$, such that $n \le (1+\ve)N$. This result
was developed in several directions. In particular, Erd\H{o}s
(see \cite{erd1} and \cite{erd2}) found a necessary and
sufficient condition on the system of nodes, for this type of
simultaneous approximation and interpolation to be valid.

We generalize the results of Bernstein and Erd\H{o}s in the following Theorem.
In order to accomplish this, we specify the choice of points $z_1,\ldots,z_N$
in an optimal fashion from the point of view of interpolation theory. Namely,
we require that the discrete measure
$$
\mu_N=\frac{1}{N}\sum_{j=1}^N\de_{z_j},
$$
where $\de_z$ denotes the unit mass placed at $z$, is close to the equilibrium measure
for $E$ (for details, see \cite{saftot1}). Fekete points (see \cite{pom1},
\cite{saftot1}) are  natural candidates for this purpose.

A Jordan curve is called quasiconformal if it is an image of the unit circle
under a quasiconformal homeomorphism of the complex plane onto itself, with
infinity as a fixed point (see \cite{lehvir} for details).

\begin{th}
\label{th2}
Let $E$ be a closed Jordan domain bounded by a quasiconformal curve
$L$. Let $f,\, r,\, k$ be as in Theorem \ref{th1} and let $z_1,\ldots,z_N\in E$
be the points of an $N$-th Fekete point set of $E$. Then for any $\ve>0$ there
exists a polynomial $p_n\in\bP_n,\, n\le (1+\ve)N,$ satisfying conditions
(\ref{1.6}) and (\ref{1.7}). Moreover, if (\ref{1.77}) holds then in
addition to (\ref{1.6}) and (\ref{1.7}) we have (\ref{1.8}), and the
constants $c_1,c_3,c_4$ and $\alpha$ are independent of $N$.
\end{th}

\absatz{Auxiliary results}

In this section, we give some results from \cite{and2}-\cite{and3}, \cite{bel},
which are needed for the proofs of the above theorems and which  characterize
the properties of the mappings $\Phi$ and $\Psi$ in the case $E\in H^*$. For
$a>0$ and $b>0$, we will use the expression $a\ole b$ (order inequality) if
$a\le cb$. The expression $a\asymp b$ means that $a\ole b$ and $b\ole a$
simultaneously. The distance $\rho_\de(z) $ to the level lines of $\Phi$ is,
for any $z\in L$, a normal majorant (in the terminology of \cite{tam}), i.e.,
\beq
\label{2.11}
\rho_{2\de}(z)\ole \rho_\de(z)\quad (\de>0).
\eeq
Let $z,\z\in L,\, \de>0$. The condition $|z-\z|\ole \rho_\de(z)$ yields
\beq
\label{2.22}
\rho_\de(\z)\asymp\rho_\de(z).
\eeq
If $L$ is a quasiconformal curve, $z\in L,\, \z\in\Om$ and if
$|z-\z|\ge\rho_\de(z)$, then the inequality
\beq
\label{2.33}
\frac{\rho_\de(z)}{|z-\z|}\ole\left(\frac{\de}{|\Phi(z)-\Phi(\z)|}\right)^\alpha
\eeq
holds with some $\alpha=\alpha(E)$.

One of the fundamental problems that, as a rule, is encountered in the construction
of approximations by polynomials, is the problem of approximating the Cauchy kernel
$1/(\z-z),\, z\in E,\, \z\in\ov{\Om}$, by polynomial kernels of the form
\beq
\label{2.44}
K_n(\z,z)=\sum_{j=0}^na_j(\z)\, z^j.
\eeq
The most general  kernels of such type, the functions $K_{r,m,k,n}(\z,z)$,
were introduced by Dzjadyk (see \cite[Chapter 9]{dzj} or
\cite[Chapter 3]{andbeldzj}).  Taking them as a basis
for our discussion, we can establish the following
result (cf. \cite[Lemma 9]{and1}).

\begin{lem}
\label{lem2.1}
Let $E\in H^*$, and let $m,\, r\in\N$.
Then for any $n\in\N$ there exists a polynomial kernel
of the form (\ref{2.44}) such that the following relations hold for
$l=0,\ldots,r$,
$z\in L$
and $\z\in\ov{\Om}$ with $d(\z,E)\le 3$:
$$
\left|\frac{\partial^l}{\partial z^l}\left(
\frac{1}{\z-z}-K_n(\z,z)\right)\right|\le\frac{c_1}{|\z-z|^{l+1}}\left(\frac{
\rho_{1/n}(z)}{|\z-z|+\rho_{1/n}(z)}\right)^m,
$$
\beq
\label{2.66}
\left|\frac{\partial^l}{\partial z^l}K_n(\z,z)\right|
\le\frac{c_2}{(|\z-z|+\rho_{1/n}(z))^{l+1}}\, ,
\eeq
where $c_j=c_j(m,r,E),\, j=1,2.$
\end{lem}
In order to improve the
approximation properties of the polynomial kernel $K_n(\z,z)$
inside of $E$, we use an idea from \cite[Theorem 2]{shi2},
completing it by the following geometrical fact. Let
$$
d(\z,B):=\mb{dist}(\z,B)=\inf\{|\z-z|: z \in B\}\quad (\z\in\C,\,B\subset\C).
$$

\begin{lem}
\label{lem2.2}
Let $E\in H^*,\, E^0\not=\emptyset$. For any $\z\in\ov{\Om}$ with
$d(\z,L)\le 3$, there exists a Jordan domain $G_\z$ with the following properties:

(i) $\z\in\partial G_\z,\, E\subset\ov{G_\z}$;

(ii) $\mb{\em diam}\,G_\z\le c\, ;$

(iii) $\partial G_\z$ is $K$-quasiconformal.

Here, the constants $c>\mb{\em diam}\,E$ and $K\ge 1$ are independent of $\z$.
\end{lem}

{\bf Proof.} If $\z\in\Om$ we set $\cal Z:=\z$; if $\z\in L$ we denote
by ${\cal Z}\in\tilde L$ the prime end whose impression coincides with $\z$
(or any of such prime ends). Let
$$
\Gamma_\z:=\{ \xi\in\Om:\, \mb{arg}\, \Phi(\xi)=\mb{arg}\, \Phi(\cal Z)\}.
$$
By virtue of \cite[Lemma 1 and Lemma 2]{and4},
\beq
\label{2.1}
d(z,L)\oge |z-\z|\quad (z\in\Gamma_\z),
\eeq
and for any $z_1,z_2\in\Gamma_\z$ the length of the part of $\Gamma_\z$
between these points satisfies
\beq
\label{2.2}
|\Gamma_\z(z_1,z_2)|\ole |z_1-z_2|.
\eeq
A result of Rickman \cite{ric} (see also \cite[p. 144]{andbeldzj})
together with (\ref{2.2}) imply that $\Gamma_\z$ is $K_1$-quasiconformal with some
$K_1\ge 1$ independent of $\z$, i.e., there exists a $K_1$-quasiconformal
mapping $F:\OC\to\OC$ such that
$$
F(\z)=0,\,\, F(\infty)=\infty,\,\, F(\Gamma_\z)=\{ w:\, w>0\}.
$$
We can assume that $|F(z_0)|=1$ for a fixed $z_0\in E^0$. We recall the following
well-known property of quasiconformal automorphisms of the complex plane
(see, for example, \cite[p. 98]{andbeldzj}): If $|\xi_1-\xi_2|
\ole |\xi_1-\xi_3|$ then
\beq
\label{2.3}
|F(\xi_1)-F(\xi_2)| \ole |F(\xi_1)-F(\xi_3)|
\eeq
and vice versa.

According to (\ref{2.1}) and (\ref{2.3}) there are constants
$c_1$ and $c_2$ such that
$$
F(E)\subset G'_\z:=\{ w=re^{i\theta}:\, 0\le r<c_1,\,
c_2<|\theta|\le\pi\}.
$$
By the Ahlfors criterion (see \cite{ahl}, \cite[p. 100]{lehvir}),
$\partial G'_\z$ is $K_2$-quasiconformal with $K_2=K_2(c_1,c_2)\ge 1$.
Therefore, by (\ref{2.3}) the domain $G_\z:=F^{-1}(G_\z')$
satisfies the conditions (i)-(iii) with $K=K_1\, K_2$.

\hfill\rbox

Let $E,\,\z$ and $G_\z$ be as in Lemma \ref{lem2.2} and let $z_0\in E^0$ be
fixed. Consider the conformal mapping $\Phi_\z:\,\OC\setminus \ov{G_\z}\to
\Delta$ normalized as in (\ref{1.55}), and the conformal mapping $\phi_\z:\,
G_\z\to \{ w:\, |w-\frac{1}{2}|<\frac{1}{2}\}$ normalized by the conditions
$$
\phi_\z(z_0)=\frac{1}{2},\,\, \phi_\z(\z)=1.
$$
Next, we use results from the theory of local distortion,
under conformal mappings of an arbitrary simply connected domain onto
a canonical one, developed by Belyi \cite{bel} (see also
\cite{andbeldzj}).

Lemma \ref{lem2.2} as well as \cite[Theorem 1 and Theorem 6]{bel}
imply that the functions $\Phi^{-1}_\z$ and $\phi_\z$ satisfy
a H\"older condition
(with constants independent of $\z$). Therefore, by \cite[Theorem 4]{bel}
for any $M\in\N$ there exists a polynomial $t_M(\z,z)\in\bP_M$
(in $z$) such that
$$
||\phi_\z-t_{M}(\z,\cdot)||_{\ov{G_\z}}\le\frac{c_1}{M^\beta}
$$
with some $c_1$ and $\beta$ independent of $\z$. We can assume that $t_M(\z,\z)=1$.

Now for $n\in\N$, we set
$$
M:=\left[\frac{n^{1/(1+\beta)}}{2}\right],\,\,
N:=[n^{\beta/(1+\beta)} ]
$$
(here $[x]$ denotes the Gauss bracket of $x$, the largest integer not exceeding $x$)
and we note that, for the polynomial
$$
u_{n/2}(\z,z):=t_M^N(\z,z),
$$
the inequality
\beq
\label{2.5}
||u_{n/2}(\z,\cdot)||_E\le\left( 1+\frac{c_1}{M^\beta}\right)^N\ole 1
\eeq
holds, as well as for any compact set $K\subset E^0$ and $\alpha:=
\beta/(1+\beta)$,
\beq
\label{2.6}
||u_{n/2}(\z,\cdot)||_K\le (1-c_2)^N\le e^{-cn^\alpha},
\eeq
where the constants $c_2<1$ and $c$ are independent of $\z$.

Hence, the function defined by
$$
T_n(\z,z):=\frac{1-u_{n/2}(\z,z)}{\z-z}+
u_{n/2}(\z,z)\, K_{[n/2]}(\z,z),
$$
where $K_{[n/2]}(\z,z)$ is the  polynomial kernel from Lemma
\ref{lem2.1}, is a polynomial (in $z$) of degree at most $n$. According to Lemma
\ref{lem2.1}, (\ref{2.5}) and (\ref{2.6}), it satisfies for $\z\in\ov{\Om},\,
d(\z,L)\le 3$, arbitrary but fixed $m\in\N$ and each compact set $K\subset E^0$
the following conditions:
$$
\left|\frac{1}{\z-z}-T_n(\z,z)\right|=|u_{n/2}(\z,z)|
\left|\frac{1}{\z-z}-K_{[n/2]}(\z,z)\right|
$$
\beq
\label{2.7}
\ole\left\{\begin{array}{ll}
\displaystyle
\frac{1}{|\z-z|}\left(\frac{\rho_{1/n}(z)}{|\z-z|+\rho_{1/n}(z)}
\right)^m, & \mb{ if } z\in L,\\[2ex]
e^{-cn^\alpha}, &\mb{ if } z\in K.
\end{array}\right.
\eeq
In addition,
\beq
\label{2.8}
|T_n(\z,z)|\ole\frac{1}{|\z-z|}\quad
(z\in E,\, \z\in\ov{\Om},\, d(\z,L)\le 3).
\eeq
We will also need the continuous extension of an arbitrary
function $F\in A(E)$ into the complex plane which preserves the smoothness
properties of $F$. The corresponding construction, proposed
by Dyn'kin \cite{dyn}, \cite{dyn1}, is based on the Whitney partition of unity (see
\cite{ste}) and local properties of the $k$-th modulus of continuity of $F$.
A slight modification of the reasoning in \cite{dyn}, \cite{dyn1} and
\cite{ste} gives the following result (cf. \cite[pp. 13-15]{andbeldzj}).
\begin{lem}
\label{lem2.3}
Let $E\in H^*$. Any $F\in A(E)$ can be continuously
extended to the complex plane (we preserve the notation $F$ for the extension)
such that:

(i) $F(z)=0$ for $z$ with $d(z,E)\ge 3$, i.e., $F$ has  compact support;

(ii) for $z\in\C\setminus E$,
$$
\left|\frac{\partial F(z)}{\partial \ov z}\right|
\le c_1\, \frac{\om_{F,k,z^*,E}(
23\, d(z,E))}{d(z,E)},
$$
where $z^*\in E$ is an arbitrary point among those ones which are
closest to $z$,
$c_1=c_1(k,\mb{\em diam}E)$;

(iii) if $\z\in E,\, z\in \C,\, |z-\z|<\de,\,
0<\de<\frac{1}{2}${\em diam}$E$, then
$$
|F(z)-P_{F,k,\z,E,\de}(z)|\le c_2\, \om_{F,k,\z,E}(25\,\de),
$$
where $P_{F,k,\z,E,\de}(z)\in\bP_{k-1}$ is the (unique) polynomial
such that
$$
||F-P_{F,k,\z,E,\de}||_{E\cap D(\z,\de)} = \om_{F,k,\z,E}(\de),
$$
and $c_2=c_2(k)$;

(iv) if $F$ satisfies a Lipschitz condition on $E$, i.e.,
$$
|F(z)-F(\z)|\le c\, |z-\z|\quad (z,\z\in E),
$$
then the extension satisfies the same condition for $z,\z\in\C$, with
$c_3=c_3(c,\mb{\em diam}E,k)$ instead of $c$.
\end{lem}

\absatz{Proof of Theorem \ref{th1}}

We fix a point $z_0\in E$ and consider a primitive of $f$:
\beq
\label{3d}
F(\z):=\int\limits_{\ga(z_0,\z)}f(\xi)\, d\xi\quad (\z\in E),
\eeq
where $\ga(z_0,\z)\subset E$ is an arbitrary rectifiable arc joining $z_0$
and $\z$.

On writing for $z\in L,\, \z\in E$ with $|\z-z|\le\de$,
\begin{eqnarray*}
F(\z)&=&F(z)+ \int\limits_{\ga(z,\z)}f(\xi)\, d\xi
\\
&=&\nu_\de(\z,z)+\int\limits_{\ga(z,\z)}\left( f(\xi)-P_{f,k,z,E,c\de}
(\xi)\right) d\xi,
\end{eqnarray*}
where $c\ge 1$ is the constant from (\ref{1.3}),
we obtain
$$
\om_{F,k+1,z,E}(\de)\le||F-\nu_\de(\cdot,z)||_{E\cap D(z,\de)}\ole\de\,
\om(\de),
$$
where $\om(\de):=\om_{f,k,E}(\de)$.
Using Lemma \ref{lem2.3}, we can extend $F$ continuously to $\C$, so that
$F$ has  compact support and satisfies
\beq
\label{3.1}
\left|\frac{\partial F(\z)}{\partial\ov{\z}}\right|\ole\om(d(\z,L)),
\eeq
for $\z\in\Om^*:=\{\z\in\ov{\Om}: \, d(\z,L)\le 3\}.$
Moreover, for $z\in L,\, \z\in\C$ with $|z-\z|\le\de<\frac{1}{2}$diam$E$,
we have
\beq
\label{3.2}
|F(\z)-\nu_\de(\z,z)|\ole\de\, \om(\de).
\eeq
Indeed, since for $\z\in E\cap D(z,\de)$,
\begin{eqnarray*}
&&|\nu_\de(\z,z)-P_{F,k+1,z,E,\de}(\z)|\\
&\le& |F(\z)-\nu_\de(\z,z)|+|F(\z)-P_{F,k+1,z,E,\de}(\z)|\ole\de\,\om(\de),
\end{eqnarray*}
we have by the Bernstein-Walsh lemma \cite[p. 77]{wal}
$$
||\nu_\de(\cdot,z)-P_{F,k+1,z,E,\de}||_{D(z,\de)}\ole
\de\,\om(\de).
$$
Hence (\ref{3.2}) follows from the last inequality and assertion (iii)
of Lemma \ref{lem2.3}.

Next, we consider the most complicated case, that is, $E^0\not= \emptyset$ and
(\ref{1.77}) holds. We introduce the polynomial kernel
$Q_{n/2}(\z,z):=T_{[n/2]}(\z,z)$, which  by (\ref{2.7})
and (\ref{2.8}) satisfies
\beq
\label{3.3}
\left|\left| \frac{1}{\z-\cdot}-Q_{n/2}(\z,\cdot)\right|\right|_K
\ole e^{-cn^\alpha} \quad (\z\in\Om^*)
\eeq
on each compact set $K\subset E^0$, and
\beq
\label{3.33}
\left| \frac{1}{\z-z}-Q_{n/2}(\z,z)\right| \ole
\frac{1}{|\z-z|}\left(\frac{\rho_{1/n}(z)}{|\z-z|+\rho_{1/n}(z)}
\right)^k\quad (z\in L),
\eeq
\beq
\label{3.34}
|Q_{n/2}(\z,z)|\ole\frac{1}{|\z-z|}\quad (z\in E).
\eeq
Further, we consider the polynomial
$$
t_n(z)=-\frac{1}{\pi}\int\limits_{\Om^*}\frac{\partial F(\z)}{\partial\ov{\z}}
\,Q_{n/2}^2(\z,z)\, dm(\z)\quad (z\in E),
$$
where $dm(\z)$ means integration with respect to the
two-dimensional Lebesgue measure
(area). Let $z\in L,\, D:=D(z,\rho),\,\sigma:=\partial D,\,
\rho:=\rho_{1/n}(z)$. According to assertion (iv) of Lemma \ref{lem2.3},
$F$ is an ACL-function (absolutely continuous on
lines parallel to the coordinate
axes) in $\C$. Hence Green's formula can be applied here (see \cite{lehvir})
to obtain
\begin{eqnarray}
f(z)-t_n(z)&=&
\frac{1}{\pi}\int\limits_{\Om^*\setminus D}
\frac{\partial F(\z)}{\partial\ov{\z}}
\,\left( Q^2_{n/2}(\z,z)-\frac{1}{(\z-z)^2}\right) dm(\z)
\nonumber
\\
&+&
\frac{1}{\pi}\int\limits_{D}
\frac{\partial F(\z)}{\partial\ov{\z}}
\, Q^2_{n/2}(\z,z)\, dm(\z)
\nonumber
\\
&+&
f(z)-\frac{1}{2\pi i}\int\limits_\sigma\frac{F(\z)}{(\z-z)^2}\,
d\z
\nonumber
\\
\label{3.4}
&=&U_1(z)+U_2(z)+U_3(z).
\end{eqnarray}
The first two integrals in (\ref{3.4}) can be estimated
in an appropriate way by passing to  polar coordinates and using
(\ref{1.5}), (\ref{1.77}), (\ref{3.1}),
(\ref{3.33})
as well as (\ref{3.34}):
\beq
\label{3.5}
|U_1(z)|\ole\int\limits_\rho^c\om(t)
\frac{\rho^{k+1}}{t^{k+2}}\, dt\ole \om(\rho)\,\rho\int\limits_\rho^c
\frac{dt}{t^2}
\ole\om(\rho),
\eeq
\beq
\label{3.6}
|U_2(z)|\ole\int\limits^\rho_0\frac{\om(t)}{t}dt\ole\om(\rho).
\eeq

In order to estimate the third term in (\ref{3.4}), we note that
$$
|f(z)-(\nu_{\rho})'_\z (z,z)|=|f(z)-P_{f,k,z,E,c\rho}(z)|\le\om(c\rho)\ole
\om(\rho),
$$
so that by (\ref{3.2}):
\beq
\label{3.7}
|U_3(z)|\le
|f(z)-(\nu_{\rho})'_\z (z,z)|+\frac{1}{2\pi}\left|\int\limits_\sigma\frac{F(\z)-
\nu_\rho(\z,z)}{(\z-z)^2}d\z\right|\ole\om(\rho).
\eeq
Comparing (\ref{3.4})-(\ref{3.7}), we obtain that
\beq
\label{3.8}
|f(z)-t_n(z)|\ole\om(\rho_{1/n}(z))\quad (z\in L).
\eeq
The estimate
\beq
\label{3.9}
||f-t_n||_K\le e^{-cn^\alpha},
\eeq
for any compact set $K\subset E^0$, follows immediately from (\ref{3.1})
and (\ref{3.3}) by a straight-forward modification of the above reasoning.

To satisfy the interpolation condition (\ref{1.7}), we argue as follows. Let $n>2N$.
We consider the polynomials
$$
V_{n/2+1}(\z,z):=
\left\{\begin{array}{ll}
1-(\z-z)Q_{n/2}(\z,z),&\mb{ if }\z\in L,\, z\in E,\\[2ex]
1,&\mb{ if } \z\in E^0,\, z\in E,
\end{array}\right.
$$
and
$$
u_n(z):=\sum_{j=1}^N\frac{q(z)}{q'(z_j)(z-z_j)}(f(z_j)-t_n(z_j))\,
V_{n/2+1}(z_j,z).
$$
By (\ref{3.3}), (\ref{3.33}), (\ref{3.8}) and (\ref{3.9}),
$$
|u_n(z)|\ole
\left\{\begin{array}{ll}
\displaystyle
\sum_j'\om(\rho_{1/n}(z_j))
\left(\frac{\rho_{1/n}(z)}{|z-z_j|+\rho_{1/n}(z)}\right)^k,&\mb{ if }
z\in L,\\[2ex]
e^{-cn^\alpha},&\mb{ if }z\in K,
\end{array}\right.
$$
where $\sum_j'$ means the sum in all $j$
with $z_j\in L$. To show that
$$
p_n(z):=t_n(z)+u_n(z)
$$
satisfies (\ref{1.6}), (\ref{1.7}) and
(\ref{1.8}), it is sufficient to prove that the
inequality
\beq
\label{3.12}
\om(\rho_{1/n}(\z))\left(\frac{\rho_{1/n}(z)}{|z-\z|+\rho_{1/n}(z)}
\right)^k\ole\om(\rho_{1/n}(z))
\eeq
holds for any $z,\z\in L$.

This relation is trivial if $|\z-z|\le \rho_{1/n}(\z)$ (cf.
(\ref{2.22})). Hence we may assume that $|\z-z|>\rho_{1/n}(\z)$. Then by
(\ref{1.5}),
$$
\om(\rho_{1/n}(\z))\left(\frac{\rho_{1/n}(z)}{|z-\z|+\rho_{1/n}(z)}
\right)^k\le
\om(|\z-z|)\left(\frac{\rho_{1/n}(z)}{|\z-z|}
\right)^k\ole \om(\rho_{1/n}(z)),
$$
which completes the proof of (\ref{3.12}).

Note that we used assumption (\ref{1.77}) only for the estimation of
$U_2(z)$ in (\ref{3.6}). If we are interested only in relations (\ref{1.6})
and (\ref{1.7}), then we need to choose in the above reasoning $Q_{n/2}(\z,z)=
K_{[n/2]}(\z,z)$, where $K_n(\z,z)$ is the polynomial kernel from Lemma \ref{lem2.1}.
Then, instead of (\ref{3.6}), we obtain by (\ref{2.66}) that
$$
|U_2(\z,z)|\ole\int\limits_0^\rho \om(t)\, \frac{t\, dt}{\rho^2}\ole
\om(\rho),
$$
and (\ref{1.77}) becomes superfluous.

\hfill\rbox

\absatz{Proof of Theorem \ref{th2d}}

Since the scheme of this proof is the same as in the proof of Theorem \ref{th1},
we describe it only briefly. We  begin with the
Taylor formula for  a primitive $F$ defined by (\ref{3d}):
$$
F(\z)=
F(z)+\sum_{j=1}^{r}\frac{f^{(j-1)}(z)}{j!}(\z-z)^j+
\frac{1}{r!}\int\limits_{\ga(z,\z)}(\z-\xi)^{r}f^{(r)}(\xi)d\xi,
$$
where $z,\z\in E$ and
 an arc $\ga(z,\z)\subset E$ joins these points and satisfies
(\ref{1.3}).
Therefore,  we have for
$z\in L,\, \z\in E$ with $|z-\z|\le\de$,
$$
F(\z)
=\ka_\de(\z,z)+\frac{1}{r!}
\int\limits_{\ga(z,\z)}(z-\xi)^r
\left( f^{(r)}(\xi)-P_{f^{(r)},k,z,E,c\de}
(\xi)\right) d\xi,
$$
where $c\ge 1$ is the constant from (\ref{1.3}) and $\ka_\de(\z,z)$
is a polynomial (in $\z$) of degree $\le k+r$. Using Lemma \ref{lem2.3},
we extend $F$ continuously, so that $F$ has  compact support and satisfies
$$
\left|\frac{\partial F(\z)}{\partial\ov{\z}}\right|\ole
d(\z,L)^r\,\om(d(\z,L))
\quad
(\z\in\Om^*:=\{\z\in\ov{\Om}:
\, d(\z,L)\le 3\}) ,
$$
$$
|F(\z)-\ka_\de(\z,z)|\ole\de^{r+1}\, \om(\de)\quad
(z\in L,\, \z\in\C,\, |\z-z|\le\de),
$$
where $\om(\de):=\om_{f^{(r)},k,z,E}(\de)$.

Next, we introduce the polynomial
$$
t_n(z)=-\frac{1}{\pi}\int\limits_{\Om^*}\frac{\partial F(\z)}{\partial\ov{\z}}
\,\frac{\partial}{\partial z} K_{n}(\z,z)\, dm(\z)\quad (z\in E),
$$
where
$K_n(\z,z)$ is the polynomial kernel from Lemma \ref{lem2.1} (with $m=2r$).

Let $l=0,\ldots,r$ and let
 $z,\, D$ as well as $\sigma$ be the same as in (\ref{3.4}).

By Green's formula, we have that
\begin{eqnarray*}
f^{(l)}(z)-t^{(l)}_n(z)&=&
\frac{1}{\pi}\int\limits_{\Om^*\setminus D}
\frac{\partial F(\z)}{\partial\ov{\z}}
\,\frac{\partial^{l+1}}{\partial z^{l+1}}
\left( K_{n}(\z,z)-\frac{1}{\z-z}\right) dm(\z)
\\
&+&
\frac{1}{\pi}\int\limits_{D}
\frac{\partial F(\z)}{\partial\ov{\z}}
\,\frac{\partial^{l+1}}{\partial z^{l+1}}
 K_{n}(\z,z)\, dm(\z)
\\
&+&
f^{(l)}(z)-\frac{1}{2\pi i}\int\limits_\sigma F(\z)
\frac{\partial^{l+1}}{\partial z^{l+1}}
\frac{1}{\z-z}\, d\z.
\end{eqnarray*}
Reasoning as in the proof of (\ref{3.8}), we obtain that
\beq
\label{1d}
|f^{(l)}(z)-t_n^{(l)}(z)|\ole\rho^{r-l}_{1/n}(z)\,
\om(\rho_{1/n}(z))\quad (z\in L).
\eeq
Further, we assume that $n>2N(r+1)$ and introduce the auxiliary polynomials
$$
V_{n/2}(\z,z):=1+\frac{(\z-z)^{r+1}}{r!}
\frac{\partial^r}{\partial z^r}K_{[n/2]}(\z,z)
$$
and
$$
u_n(z):=\sum_{j=1}^N\frac{q^{r+1}(z)}{(z-z_j)^{r+1}}\, V_{n/2}(z_j,z)\,
\sum_{s=0}^rA_{j,s}\, (z-z_j)^s,
$$
where
$$
A_{j,s}:=\left.\sum_{\nu=0}^s\frac{1}{\nu!(s-\nu)!}\left( f^{(\nu)}(z_j)-
t_n^{(\nu)}(z_j)\right)\left(\frac{\partial^{s-\nu}}{\partial z^{s-\nu}}
\frac{(z-z_j)^{r+1}}{q^{r+1}(z)}\right)\right|_{z=z_j}.
$$
According to the Hermite interpolation formula (see \cite{smileb}),
we have
$$
u_n^{(l)}(z_j)=f^{(l)}(z_j)-t_n^{(l)}(z_j)\quad (j=1,\ldots,N).
$$
Therefore the polynomial
$$
p_n:=u_n+t_n
$$
satisfies the interpolation condition
(\ref{1.7d}).

Since
$$
|A_{j,s}|\ole \rho_{1/n}^{r-s}(z_j)\, \om(\rho_{1/n}(z_j)),
$$
we obtain by Lemma \ref{lem2.1} for any $z\in L$,
\begin{eqnarray}
|u_n(z)|&\ole&
\sum_{j=1}^N\left(\frac{\rho_{1/n}(z)}{|z-z_j|+\rho_{1/n}(z)}
\right)^{2r}\sum_{s=0}^r \rho_{1/n}^{r-s}(z_j)\,
\om(\rho_{1/n}(z_j))\, |z-z_j|^s
\nonumber
\\
\label{4d}
&\ole& \rho_{1/n}^r(z)\, \om(\rho_{1/n}(z)),
\end{eqnarray}
where we used (\ref{2.22}) and the following inequality:
for $z,\,\z\in L$ with $|\z-z|\ge\rho_{1/n}(z)$,
$$
\left|\frac{\rho_{1/n}(z)}{z-\z}\right|^{2r}|z-\z|^r\,
\om(|z-\z|)\ole
\rho_{1/n}^r(z)\, \om(\rho_{1/n}(z)).
$$
By a theorem of Tamrazov \cite{tam} (see also \cite[p. 187]{andbeldzj}),
(\ref{4d}) yields
\beq
\label{2d}
|u_n^{(l)}(z)|\ole \rho_{1/n}^{r-l}(z)\, \om(\rho_{1/n}(z)).
\eeq
Combining (\ref{1d}) and (\ref{2d}), we obtain (\ref{1.6d}).

\hfill\rbox

\absatz{Proof of Theorem \ref{th2}}

We use the same scheme as in the proof of Theorem \ref{th1}. Let
(\ref{1.77}) hold. We
construct a polynomial $t_N\in\bP_N$ such that
\beq
\label{4.1}
|f(z)-t_N(z)|\ole\om(\rho_{1/N}(z))\quad (z\in L),
\eeq
where $\om(\de):=\om_{f,k,E}(\de)$, and
\beq
\label{4.2}
||f-t_N||_K\le e^{-cN^\alpha}
\eeq
for any compact set $K\subset E^0$.

Let $m:=[\ve N]$. Consider the polynomial
$$
u_{N+m}(z):=\sum_{j=1}^N\frac{q(z)}{q'(z_j)(z-z_j)}(f(z_j)-t_N(z_j))V_{m+1}(z_j,
z),
$$
where
$$
V_{m+1}(\z,z):=1-(\z-z)Q_m(\z,z)\quad (\z\in L,\, z\in E),
$$
and $Q_m(\z,z):=T_m(\z,z)$ is a polynomial of degree at most $m$ (in $z$)
satisfying the  inequalities (cf. (\ref{2.7}))
\beq
\label{4.3}
\left|\frac{1}{\z-z}-Q_m(\z,z)\right|\ole\frac{1}{|\z-z|}\left(\frac
{\rho_{1/m}(z)}{|\z-z|+\rho_{1/m}(z)}\right)^{k+l}\quad
(z,\z\in L)
\eeq
(the choice of $l=l(E)>0$ will be specified below)
and
\beq
\label{4.4}
\left|\left|\frac{1}{\z-\cdot}-Q_m(\z,\cdot)\right|\right|_K\le
e^{-cm^\alpha}\quad (\z\in L)
\eeq
on each compact set $K\subset E^0$.

Let $z\in L,\, \Phi(z)=e^{i\theta_0},\, \Phi(z_j)=e^{i\theta_j},$
$$
0\le\theta_1<\theta_2<\ldots<\theta_N<\theta_{N+1}:=\theta_1+2\pi.
$$
It is proved in \cite{andbla} that
\beq
\label{4.5}
|\theta_{j+1}-\theta_j|\asymp\frac{1}{N}\quad (j=1,\ldots,N).
\eeq
We rename the points $\{ e^{i\theta_j}\}_1^N$ by
$\{ e^{i\theta'_j}\}_1^\mu$, $\{ e^{i\theta''_j}\}_1^\nu$ and
$\{ e^{i\theta'''_j}\}_1^{N-\mu-\nu}$ in such a way that
$$
|\theta_0-\theta_j'|\le\frac{1}{m}\quad (j=1,\ldots,\mu),
$$
and $\theta_j=\theta_j'',\theta_j'''$ satisfy
$$
|\theta_0-\theta_j|>\frac{1}{m},\quad (\theta_j\not\in
\{\theta'_1,\ldots,\theta'_\mu\}),
$$
$$
\theta_0<\theta_1'' <\theta_2'' <\ldots<\theta_\nu'' \le\pi+\theta_0,
$$
$$
\theta_0-\pi<\theta_{N-\mu-\nu}'''<\ldots<\theta_{1}'''<\theta_0.
$$
Equation (\ref{4.5}) implies that
$$
\mu\asymp\frac{1}{\ve},\,\, \nu\asymp N-\mu-\nu\asymp N.
$$
Furthermore, for the function
$$
h(\theta,\theta_0):=
(f(\Psi(e^{i\theta}))-t_N(\Psi(e^{i\theta})))
V_{m+1}(\Psi(e^{i\theta}),\Psi(e^{i\theta_0}))
$$
we have by (\ref{1.5}), (\ref{2.22}), (\ref{4.1}) and (\ref{4.3}),
\beq
\label{4.6}
|h(\theta_j',\theta_0)|\ole\om(\rho),
\eeq
\beq
\label{4.7}
|h(\theta''_j,\theta_0)|\ole\om(|z-z''_j|)
\left(\frac{\rho}{|z-z''_j|}\right)^{k+l}\ole
\om(\rho)\left(\frac{\rho}{|z-z''_j|}\right)^{l},
\eeq
\beq
\label{4.8}
|h(\theta'''_j,\theta_0)|\ole\om(|z-z'''_j|)
\left(\frac{\rho}{|z-z'''_j|}\right)^{k+l}\ole
\om(\rho)\left(\frac{\rho}{|z-z'''_j|}\right)^{l},
\eeq
where $\rho:=\rho_{1/m}(z),\, z_j'':=\Psi (e^{i\theta_j''}),\,
z_j''':=\Psi (e^{i\theta_j'''})$.

It follows from (\ref{4.2}) and (\ref{4.4}) that the polynomial
\beq
\label{4.9}
p_{[(1+\ve)N]}(z):= t_N(z)+u_{N+m}(z)
\eeq
satisfies (\ref{1.7}) and (\ref{1.8}).

We choose $l$ so that
$$
\left|\frac{\rho}{\z-z}\right|^l\ole\left(\frac{1}{m|\Phi(\z)-\Phi(z)|}
\right)^2,
$$
for $\z\in L$ with $|\z-z|>\rho$ (cf. (\ref{2.33})).

Since
\begin{eqnarray*}
|u_{N+m}(z)|&\le&
\sum_{j=1}^\mu|h(\theta_j',\theta_0)|
+\sum_{j=1}^\nu|h(\theta_j'',\theta_0)|+
\sum_{j=1}^{N-\mu-\nu}|h(\theta_j''',\theta_0)|\\
&\ole&\om(\rho)\left( 1+\sum_{j=1}^N\frac{1}{j^2}\right)\ole\om(\rho), \quad z\in L,
\end{eqnarray*}
by (\ref{4.6})-(\ref{4.8}), we obtain the desired
inequality (\ref{1.6}) by (\ref{2.11}) and
(\ref{4.1}), for $p_{[(1+\ve)N]}$ given by (\ref{4.9}).

Taking in the above argument $Q_m(\z,z):=K_m(\z,z)$, we obtain equations
(\ref{1.6}) and (\ref{1.7}) even without assumption (\ref{1.77}).

\hfill\rbox

GSF-Forschungszentrum, Institut f\"ur Biomathematik und
Biometrie, Ingolst\"adter Landstr. 1, D-85764 Neuherberg, Germany;\\
e-mail: mgk002@eo-dec-mathsrv.ku-eichstaett.de

Department of Mathematics, 401 Mathematical Sciences, Oklahoma State University,
Stillwater, OK 74078-1058, U.S.A.;\\
e-mail: igor@math.okstate.edu

Institute for Computational Mathematics, Kent State University, Kent, OH 44242,
U.S.A.;\\
e-mail: varga@mcs.kent.edu

\end{document}